\documentclass[a4paper,12pt,reqno]{amsart}
\usepackage{amsmath,amsthm,amssymb}
\usepackage{graphicx, color}
\usepackage{url}
\usepackage[numbers,sort&compress]{natbib}
\nonstopmode \numberwithin{equation}{section}

\allowdisplaybreaks

\allowdisplaybreaks

\begin{document}

\title[{Analytic computations of digamma functions }]{Analytic computations of digamma function\\ using some new identities}

\author[M.I. Qureshi and M. Shadab]{Mohammad Idris Qureshi and Mohd Shadab$^{*}$}

\address{Mohammad Idris Qureshi:   Department of Applied Sciences and Humanities,
 Faculty of Engineering and Technology,
 Jamia Millia Islamia (A Central University),
 New Delhi 110025, India}
\email{miqureshi\_delhi@yahoo.co.in}

\address{Mohd Shadab: Department of Applied Sciences and Humanities,
 Faculty of Engineering and Technology,
 Jamia Millia Islamia (A Central University),
 New Delhi 110025, India}
\email{shadabmohd786@gmail.com}



\thanks{*Corresponding author}

\begin{abstract}
Motivated by rigorous development in the theory of digamma functions, we have first derived some new identities for the digamma function, and then computed the values of digamma function for the fractional orders using these identities conveniently.
 \end{abstract}

\maketitle

\section{Introduction and preliminaries}\label{Intro}
A natural property of digamma (Psi function) function to be used as application in the theory of beta distributions-probability models for the domain [0,1]. It is mainly used in the theory of special functions with a wide range of the applications. Digamma functions are directly connected with many special functions such as Riemann's zeta function and Clausen's function etc.\\

 Many authors have contributed to develop the theory of polygamma function with respect to properties \cite{Qi, Gauss, Jensen, Lehmer, Mahler}, inequalities \cite{Alzer, Batir, Clark}, monotonicity \cite{Qi1, Qi2, Qi3, Qi4}, series \cite{Bor, de, Lewin, Wu, Gosper, Grossman}, and fractional calculus \cite{Adam, Al-saqabi, Sri1, Sri2}.\\

The Gamma function, $\Gamma{(z)}$, was introduced by Leonard Euler as a generalization of the factorial function on the sets, $\mathbb{R}$ of all real numbers, and $\mathbb{C}$ of all complex numbers. It (or, Euler's integral of second kind) is defined by

\begin{eqnarray}\label{eq(1.1)}
\Gamma{(z)}&=&\int_{0}^{\infty}\exp{(-t)}t^{z-1}dt,\qquad \Re(z)>0\nonumber\\
&=&\lim_{n\rightarrow\infty}\int_{0}^{n}{\left(1-\frac{t}{n}\right)^n}t^{z-1}dt.
\end{eqnarray}
\vskip.2cm
In 1856, Karl Weierstrass gave a novel definition of gamma function
\begin{eqnarray}\label{eq(1.2)}
\frac{1}{\Gamma{(z)}}&=&z \exp{(\gamma z)} \prod_{n=1}^{\infty}\left[\left(1+\frac{z}{n}\right)\exp{\left(-\frac{z}{n}\right)}\right],
\end{eqnarray}
where $\gamma= 0.577215664901532860606512090082402431042\dots$ is called Euler-Mascheroni constant, and
$\frac{1}{\Gamma{(z)}}$ is an entire function of $z$, and
\begin{eqnarray}
\gamma&=&\lim_{n\rightarrow\infty}\left(1+\frac{1}{2}+\frac{1}{3}+.....+\frac{1}{n}-\ell n{(n)}\right)\nonumber.
\end{eqnarray}

The function

\begin{eqnarray}
\psi(z)=\frac{d}{dz}\{\ell n \,{\Gamma(z)}\}=\frac{\Gamma ^{\prime}(z)}{\Gamma(z)},
\end{eqnarray}
or, equivalently
\begin{eqnarray}
\ell n \,{\Gamma(z)}=\int_{1}^{z}\psi(\zeta)d\zeta,
\end{eqnarray}

is the logarithmic derivative of the gamma function (Psi function or digamma function).\\

$\psi^{(i)}(z)$ for $i\in \mathbb{N}$ are called the polygamma functions, and $\psi$ has the presentation as
\begin{eqnarray}
\psi(z)=\frac{\Gamma ^{\prime}(z)}{\Gamma(z)}
=-\gamma+\int_{0}^{\infty}\frac{e^{-t}-e^{-zt}}{1-e^{-t}}dt  \quad (\gamma= \text{ Euler's constant}).
\end{eqnarray}

The Psi function has following series representation
\begin{eqnarray}
\psi(z)=-\gamma-\frac{1}{z}+\sum_{n=1}^{\infty}\frac{z}{n(z+n)}, \quad z \ne 0,-1,-2,-3,\dots
\end{eqnarray}

The {\it generalized hypergeometric function} ${_p}F_q$, is defined by

\begin{eqnarray}\label{eq.(1.11)}
{_p}F_{q} \left[ \begin{array}{r}(a_p); \\ (b_q); \end{array} z \right] = \sum_{m=0}^{\infty} \frac{[(a_p)]_m}{[(b_q)]_m} \frac{z^m}{m!} ,
\end{eqnarray}

\begin{itemize}
  \item p and q are positive integers or zero,
  \item $z$ is a complex variable,
  \item $(a_p)$ designates the set ${a_1,a_2, . . . , a_p}$,
  \item the numerator parameters $a_1, . . . , a_p \in \mathbb{C}$ and the denominator parameters $b_1, . . . , b_q \in \mathbb{C}\setminus {\mathbb{Z}^{-}_{0}}, $
  \item $[(a_r)]_k = \displaystyle\prod_{i=1}^{r} (a_i)_k$. By convention, a product over the empty set is 1,
  \item $(a)_k$ is the Pochhammer's symbol.
\end{itemize}

The widely used Pochhammer symbol $(\lambda)_{\nu}$ ~$(\lambda, \nu \in\mathbb{C})$ is defined by
\begin{equation}
\left(\lambda\right)_{\nu}:=\frac{\Gamma\left(\lambda+\nu\right)}{\Gamma\left(\lambda\right)}=\begin{cases}
\begin{array}{c}
1\\
~\\
\lambda\left(\lambda+1\right)\ldots\left(\lambda+n-1\right)
\end{array} & \begin{array}{c}
\left(\nu=0;\lambda\in\mathbb{C}\setminus\left\{ 0\right\} \right)\\
~\\
\left(\nu=n\in\mathbb{N};\lambda\in\mathbb{C}\right),
\end{array}\end{cases}
\end{equation}
it being understood $conventionally$ that $\left(0\right)_{0}=1$, and assumed $tacitly$ that the $\Gamma$ quotient exists.
\vskip.2cm
Thus, if a numerator parameter is a negative integer or zero, the $_pF_q$ series terminates, and then we are led to a generalized hypergeometric polynomial.
\vskip.2cm

In 1813, Gauss \cite{Gauss} (see also, Jensen \cite[p.146, eq.(32)]{Jensen};  \cite[p.19, (1.7.3) eq.(29)]{Erde}; B\"{o}hmer \cite[p.77]{Boh} ) discovered an interesting formula for digamma (Psi) function as follows

\begin{eqnarray}\label{eq(Gauss)}
\psi(p/q) = -\gamma - \ell n\,{(q)}-\frac{\pi}{2}\cot{\left(\frac{\pi p}{q}\right)}+\sum_{j=1}^{[\frac{q}{2}]}{'}\left\{\cos\left({\frac{2\pi j p}{q}}\right)\ell n\, {\left(2-2\cos{\frac{2\pi j}{q}}\right)}\right\},
 \end{eqnarray}

where $1 \le p < q$ and $p, q$ are positive integers, and accent(prime) to right of the summation sign indicates the term corresponding to (last term) $j = \frac{q}{2}$ (when $q$ is positive even integer) should be divided by 2.\\

A different form of Gauss formula is also given in N. Nielsen \cite[p. 22, an equation between equations (7) and (8)]{Nielsen} as follows

\begin{eqnarray}\label{eq(Gauss1)}
\psi(p/q) = -\gamma - \ell n\,{(q)}-\frac{\pi}{2}\cot{\left(\frac{\pi p}{q}\right)}+\sum_{j=1}^{q-1}\left\{\cos\left({\frac{2\pi p j }{q}}\right)\ell n\, \left(2\sin\left(\frac{\pi j}{q}\right)\right)\right\},
 \end{eqnarray}

where $1 \le p < q$ and $p, q$ are positive integers.

\vskip.2cm

Afterwards, in 2007, a simplified treatment of the above formula was made by Murty and Saradha \cite[p. 300, after eq.(4)]{Murty} (see also, Lehmer \cite[p. 135, after eq.(20)]{Lehmer}) ) such that

\begin{eqnarray}\label{eq(Murty)}
\psi(p/q)=-\gamma-\ell n\,{(2q)}-\frac{\pi}{2}\cot{\left(\frac{\pi p}{q}\right)}+2\sum_{j=1}^{[\frac{q}{2}]}\left\{\cos{\left(\frac{2\pi pj}{q}\right)}\ell n\,{ \sin{\left(\frac{\pi j}{q} \right)}}\right\},
\end{eqnarray}
where $p = 1, 2, 3, \dots ,(q-1) ,\,  q = 2, 3, 4, \dots ; (p,q)= 1$.\\

Also, we have verified the results \eqref{eq(Gauss)}, \eqref{eq(Gauss1)} and \eqref{eq(Murty)} by taking different values of $p$ and $q$.

An erroneous formula of digamma function is also recorded in Gradshteyn and Ryzhik \cite[p. 904, eq 8.363(6)]{Gradshteyn} such that

\begin{eqnarray}\label{eq(Err)}
\psi(p/q) \circeq -\gamma-\ell n\,{(2q)}-\frac{\pi}{2}\cot{\left(\frac{\pi p}{q}\right)}+2\sum_{j=1}^{[\frac{q+1}{2}]-1
}\left\{\cos{\left(\frac{2\pi pj}{q}\right)}\ell n\,{ \sin{\left(\frac{\pi j}{q} \right)}}\right\},
\end{eqnarray}
where $p = 1, 2, 3, \dots ,(q-1) ,\,  q = 2, 3, 4, \dots ; (p,q)= 1$ and the symbol $\circeq$ exhibits the fact that equation \eqref{eq(Err)} does not hold true as stated.

{\bf Some important facts:}
\begin{itemize}
\item We can not compute the value of digamma function when $p > q$ or ( and ) $\frac{p}{q}$ is negative fraction using Gauss formula \cite{Gauss}.

\item We can not compute the value of digamma function when $p > q$
using Jensen formula \cite{Jensen}.

\item We can not compute the value of digamma function when $\frac{p}{q}$ is negative
using Jensen \cite{Jensen}.

\item  Murty and Saradha \cite[p. 300]{Murty} corrected a formula of Lehmer \cite[p. 135]{Lehmer} for $\psi(\frac{p}{q})$ .

\item The value of digamma function has been proved transcendental with the help of Gauss formula in \cite{Murty}.
\end{itemize}

\vskip.2cm
\section{Some new identities for digamma function}
Here, we derive some interesting  and new identities for the computation of digamma function of fractional order.
\vskip.2cm

Some functional relations for digamma function that are easily derivable from the properties of the gamma function. Indeed, from the formula,

\begin{eqnarray}\label{eq(main-2.1)}
\Gamma(z+1)=z\Gamma(z),\quad
\Gamma(z)\Gamma(1-z)=\frac{\pi}{\sin{(\pi z)}},  \quad z\ne0,\pm1,\pm2,\pm3,\dots
\end{eqnarray}
taking $\ell n$ both sides and differentiating the above equation with respect to z, we derive the some basic identities for digamma function as follows
\begin{eqnarray}\label{eq(main-2.2)}
\psi(z+1) = \psi(z)+\frac{1}{z}, \quad
\psi(1-z) = \psi(z)+\pi \cot{(\pi z)},  \quad z\ne0,\pm1,\pm2,\pm3,\dots
\end{eqnarray}

\begin{eqnarray}\label{eq(main-2.3)}
\psi(z+n) = \frac{1}{z} + \frac{1}{z+1} + \dots  +\frac{1}{z+n-1} + \psi(z).
\end{eqnarray}

On setting $z=1-z$ in equation \eqref{eq(main-2.2)}, we get
\begin{eqnarray}\label{eq(main-2.4)}
\psi(-z) =  \frac{1}{z} + \psi(1-z).
\end{eqnarray}

On comparing the values of $\psi(1-z)$ from the equations \eqref{eq(main-2.2)} and \eqref{eq(main-2.4)}, we get a new identity

\begin{equation}\label{eq(main-2.5)}
\psi(z)+\pi \cot{(\pi z)} =  \psi(-z) - \frac{1}{z}.
\end{equation}

By setting $z=\frac{p}{q}, \,\, 1\leq p < q$ in equations \eqref{eq(main-2.2)} and \eqref{eq(main-2.4)}, we get more identities, which would be used to derive our main identities to compute the values of digamma function

\begin{eqnarray}\label{eq(dig-2.6)}
\psi\left(\frac{p+q}{q}\right) =  \frac{q}{p} + \psi\left(\frac{p}{q}\right), \text{  and  }
\psi\left(\frac{-p}{q}\right) =  \frac{q}{p} + \psi\left(\frac{q-p}{q}\right), \quad 1\leq p < q.
\end{eqnarray}

\vskip.2cm

For the sake of convenient computation of digamma function, we derive some more identities, which are simple but more applicable in the computation of digamma function for $\frac{p}{q} > 1$. For this concern, we connect the Murty and Saradha's corrected formula for digamma function \eqref{eq(Murty)} with our above derived identity \eqref{eq(dig-2.6)}, and get the result as follows
\begin{eqnarray}\label{eq(main-2.7)}
\psi\left(\frac{q-p}{q}\right) = -\gamma-\ell n{(2q)}+\frac{\pi}{2}\cot{\left(\frac{\pi p}{q}\right)}+2\sum_{j=1}^{[\frac{q}{2}]}\left\{\cos{\left(\frac{2\pi pj}{q}\right)}\ell n{ \sin{\left(\frac{\pi j}{q} \right)}}\right\},
\end{eqnarray}
$(p,q)=1$; $1\leq p < q$.

\vskip.2cm

Now, we derive the identity for the computation of digamma function for negative fractions ($-\frac{p}{q}$). For this motive, we derive the identity in the similar manner as used in above identity, and get the result as follows
\begin{eqnarray}\label{eq(main-2.8)}
\psi\left(\frac{-p}{q}\right) = \frac{q}{p} - \gamma-\ell n{(2q)} - \frac{\pi}{2}\cot{\left(\frac{\pi (q-p)}{q}\right)}+2\sum_{j=1}^{[\frac{q}{2}]}\left\{\cos{\left(\frac{2\pi (q-p)j}{q}\right)}\ell n{ \sin{\left(\frac{\pi j}{q} \right)}}\right\},\nonumber\\
\end{eqnarray}
$(p,q)=1$; $1\leq p < q$.
\section{Numeric computations of digamma function}

\begin{table}[h!]
\centering

\caption{$\psi$- Function(Negative Fractional Valued)}
\vskip.1cm
\label{table-1.2}
\begin{tabular}{|c  |c  |c  |c|c|c|c|c|c|}
\hline
Ser. No. &$z=\frac{p}{q}$ & $\psi(z)= \frac{\Gamma^{'}{(z)}}{\Gamma{(z)}}$  \\ [5pt]\hline

1 &$-\frac{2}{3}$ & $ -\gamma+\frac{3}{2} - \frac{\pi \sqrt{3}}{6} - \frac{3\ell n \,3}{2}$  \\[5pt] \hline

2 &$-\frac{3}{4}$ & $-\gamma+\frac{4}{3} - \frac{\pi}{2} -3\ell n \,2 $  \\[5pt] \hline

3 &$-\frac{1}{2}$ & $-\gamma +2 - 2\ell n \,2 $   \\ [5pt] \hline

4 &$-\frac{1}{3}$ & $-\gamma +3 + \frac{\sqrt{3}\pi}{6} - \frac{3}{2}\ell n \,3 $   \\ [5pt] \hline

5 &$-\frac{1}{4}$ & $-\gamma +4 + \frac{\pi}{2} - 3\ell n \,2 $   \\ [5pt] \hline

6 &$-\frac{5}{8}$ & $-\gamma +\frac{8}{5} - \frac{(\sqrt{2} -1)\pi}{2} -4\ell n \,2 + \sqrt{2}\ell n \,(1+\sqrt{2}) $  \\[5pt] \hline

7 &$-\frac{3}{8}$ & $-\gamma+\frac{8}{3} + \frac{(\sqrt{2}-1)\pi}{2} -4\ell n \,2 + \sqrt{2}\ell n \,(1+\sqrt{2}) $  \\[5pt] \hline

8 &$-\frac{1}{8}$ & $-\gamma +8 + \frac{(1+\sqrt{2})\pi}{2} - 4\ell n \,2 -\sqrt{2}\ell n \,(1+\sqrt{2})$   \\ [5pt] \hline

9 &$-\frac{5}{6}$ & $-\gamma +\frac{6}{5} - \frac{\pi \sqrt{3}}{2} - \frac{3}{2}\ell n \,3 - 2 \ell n \,2$   \\ [5pt] \hline

10 &$-\frac{3}{2}$ & $-\gamma+\frac{8}{3} -2\ell n \,2 $  \\[5pt] \hline

11 &$-\frac{7}{3}$ & $-\gamma+\frac{117}{28} +\frac{\pi\sqrt{3}}{6}-\frac{3}{2}\ell n \,3 $  \\[5pt] \hline

\end{tabular}
\end{table}

\newpage

\begin{table}[h!]
\centering

\caption{$\psi$- Function(Positive Fractional Valued)}
\vskip.1cm
\label{table-1.2}
\begin{tabular}{|c  |c  |c  |c|c|c|c|c|c|}
\hline
Ser. No. & $z=\frac{p}{q}$ & $\psi(z)= \frac{\Gamma^{'}{(z)}}{\Gamma{(z)}}$  \\ [5pt]\hline

1 &$\frac{1}{2}$ & $-\gamma - 2\ell n \, 2 $   \\ [5pt] \hline

2 &$\frac{1}{3}$ & $-\gamma -\frac{\sqrt{3} \pi}{6}- \frac{3}{2}\ell n \,3 $   \\ [5pt] \hline

3 &$\frac{1}{4}$ & $-\gamma - \frac{\pi}{2} - 3\ell n \,2$   \\ [5pt]\hline

4 &$\frac{1}{5}$ & $-\gamma -\ell n 10 - \left(\frac{1+\sqrt{5}}{\sqrt{(10-2\sqrt{5})}} \right) \frac{\pi}{2} +  \frac{1}{2}\{ \sqrt{5} \ell n \,\left(\frac{\sqrt{5} - 1}{2}\right) - \ell n \,\frac{\sqrt{5}}{4} \}$  \\[5pt] \hline

5 &$\frac{1}{6}$ & $-\gamma -\ell n \, 12 - \frac{\pi \sqrt{3}}{2} - \ell n \, \sqrt{3} $  \\[5pt] \hline

6 &$\frac{1}{8}$ & $-\gamma - \frac{(1+\sqrt{2})\pi}{2} -4\ell n \,2 - \sqrt{2}\ell n \,(1+\sqrt{2})$   \\[5pt] \hline

7 &$\frac{1}{10}$ & $-\gamma -\ell n \, 20 - \left(\frac{\sqrt{(10+2\sqrt{5})}}{\sqrt{5} -1} \right) \frac{\pi}{2} +  \frac{1}{2}\{ \sqrt{5} \ell n \, (\sqrt{5} - 2) - \ell n \, \sqrt{5} \}$  \\[5pt] \hline

8 &$\frac{1}{12}$ & $-\gamma -\ell n \, 24 - \left(2+\sqrt{3} \right) \frac{\pi}{2} +  \{ \sqrt{3} \ell n \, (2-\sqrt{3}) - \ell n \, \sqrt{3} \}$  \\[5pt] \hline

9 &$\frac{2}{3}$ & $ -\gamma +\frac{\sqrt{3} \pi}{6} - \frac{3}{2} \ell n \, 3$   \\ [5pt] \hline

10 &$\frac{2}{5}$ & $-\gamma -\ell n 10 - \left(\frac{\sqrt{5} -1}{\sqrt{(10+2\sqrt{5})}} \right) \frac{\pi}{2} +  \frac{1}{2}\{ \sqrt{5} \ell n \, \left(\frac{\sqrt{5} + 1}{2}\right) - \ell n \, \frac {\sqrt{5}}{4} \}$  \\[5pt] \hline

11 &$\frac{3}{4}$ & $-\gamma + \frac{\pi}{2} -3\ell n \,2 $  \\[5pt] \hline

12 &$\frac{3}{5}$ & $-\gamma -\ell n 10 + \left(\frac{\sqrt{5} -1}{\sqrt{(10+2\sqrt{5})}} \right) \frac{\pi}{2} +  \frac{1}{2}\{ \sqrt{5} \ell n \, \left(\frac{\sqrt{5} + 1}{2}\right) - \ell n \, \frac {\sqrt{5}}{4} \}$  \\[5pt] \hline

13 &$\frac{3}{8}$ & $-\gamma - \frac{(\sqrt{2}-1)\pi}{2} - 4\ell n \,2 + \sqrt{2}\ell n \,(1+\sqrt{2}) $  \\[5pt] \hline

14 &$\frac{3}{10}$ & $-\gamma -\ell n \, 20 - \left(\frac{\sqrt{(10-2\sqrt{5})}}{1+\sqrt{5}} \right) \frac{\pi}{2} +  \frac{1}{2}\{ \sqrt{5} \ell n \, (2+ \sqrt{5}) - \ell n \, \sqrt{5} \}$  \\[5pt] \hline

15 &$\frac{4}{5}$ & $-\gamma -\ell n 10 + \left(\frac{1+\sqrt{5}}{\sqrt{(10-2\sqrt{5})}} \right) \frac{\pi}{2} +  \frac{1}{2}\{ \sqrt{5} \ell n \,\left(\frac{\sqrt{5} - 1}{2}\right) - \ell n \,\frac{\sqrt{5}}{4} \}$  \\[5pt] \hline

16 &$\frac{5}{6}$ & $-\gamma -\ell n \, 12 + \frac{\pi \sqrt{3}}{2} - \ell n \, \sqrt{3} $  \\[5pt] \hline

17 &$\frac{5}{8}$ & $-\gamma + \frac{(\sqrt{2}-1)\pi}{2} - 4\ell n \,2 + \sqrt{2}\ell n \,(1+\sqrt{2}) $  \\[5pt] \hline

18 &$\frac{5}{12}$ & $-\gamma -\ell n \, 24 - \left(2-\sqrt{3} \right) \frac{\pi}{2} +  \{ \sqrt{3} \ell n \, (2+\sqrt{3}) - \ell n \, \sqrt{3} \}$  \\[5pt] \hline

19 &$\frac{7}{8}$ & $-\gamma + \frac{(1+\sqrt{2})\pi}{2} - 4\ell n \,2 - \sqrt{2}\ell n \,(1+\sqrt{2}) $  \\[5pt] \hline

20 &$\frac{7}{10}$ & $-\gamma -\ell n \, 20 + \left(\frac{\sqrt{(10-2\sqrt{5})}}{1+\sqrt{5}} \right) \frac{\pi}{2} +  \frac{1}{2}\{ \sqrt{5} \ell n \, (2+ \sqrt{5}) - \ell n \, \sqrt{5} \}$  \\[5pt] \hline

21 &$\frac{7}{12}$ & $-\gamma -\ell n \, 24 + \left(2-\sqrt{3} \right) \frac{\pi}{2} +  \{ \sqrt{3} \ell n \, (2+\sqrt{3}) - \ell n \, \sqrt{3} \}$  \\[5pt] \hline

22 &$\frac{9}{10}$ & $-\gamma -\ell n \, 20 + \left(\frac{\sqrt{(10+2\sqrt{5})}}{\sqrt{5} -1} \right) \frac{\pi}{2} +  \frac{1}{2}\{ \sqrt{5} \ell n \, (\sqrt{5} - 2) - \ell n \, \sqrt{5} \}$  \\[5pt] \hline

23 &$\frac{11}{12}$ & $-\gamma -\ell n \, 24 + \left(2+\sqrt{3} \right) \frac{\pi}{2} +  \{ \sqrt{3} \ell n \, (2-\sqrt{3}) - \ell n \, \sqrt{3} \}$  \\[5pt] \hline

\end{tabular}
\end{table}

\newpage

\begin{table}[h!]
\centering

\caption{$\psi$- Function(Fractional Valued, $p>q$ )}
\vskip.1cm
\label{table-1.2}
\begin{tabular}{|c  |c  |c  |c|c|c|c|c|c|}
\hline
Ser. No. &$z=\frac{p}{q}$ & $\psi(z)= \frac{\Gamma^{'}{(z)}}{\Gamma{(z)}}$  \\ [5pt]\hline

1 &$\frac{7}{3}$ & $-\gamma+\frac{15}{4} -\frac{ \pi \sqrt{3}}{6}- \frac{3}{2}\ell n \,3 $   \\ [5pt] \hline

2 &$\frac{3}{2}$ & $ -\gamma +2 - 2 \ell n \, 2$   \\ [5pt] \hline

3 &$\frac{5}{2}$ & $-\gamma+\frac{8}{3} -2\ell n \,2 $  \\ [5pt] \hline

\end{tabular}
\end{table}

The following errata are found in a paper of Jensen \cite[p. 147]{Jensen} such that

\begin{eqnarray}\label{eq(Err1)}
\psi(3/5) \circeq -\gamma  + \frac{\pi}{2} \sqrt{\left(1-\frac{2}{\sqrt{5}} \right)}- \frac{5}{4}  \ell n \, 5 +\frac{\sqrt{5}}{4}  \ell n\, \left(\frac{3+2\sqrt{5}}{2} \right),
 \end{eqnarray}

\begin{eqnarray}\label{eq(Err2)}
\psi(4/5) \circeq -\gamma  + \frac{\pi}{2} \sqrt{\left(1+\frac{2}{\sqrt{5}} \right)}- \frac{5}{4}  \ell n \, 5 -\frac{\sqrt{5}}{4}  \ell n\, \left(\frac{3+2\sqrt{5}}{2} \right),
 \end{eqnarray}

where the symbol $\circeq$ exhibits the fact that each of the equations \eqref{eq(Err1)} and \eqref{eq(Err2)} does not hold true as stated.

The following are the corrected forms of above equations
\begin{eqnarray}
\psi(3/5) = -\gamma  + \frac{\pi}{2} \sqrt{\left(1-\frac{2}{\sqrt{5}} \right)}- \frac{5}{4}  \ell n \, 5 +\frac{\sqrt{5}}{4}  \ell n\, \left(\frac{3+\sqrt{5}}{2} \right),
 \end{eqnarray}

\begin{eqnarray}
\psi(4/5) = -\gamma  + \frac{\pi}{2} \sqrt{\left(1+\frac{2}{\sqrt{5}} \right)}- \frac{5}{4}  \ell n \, 5 -\frac{\sqrt{5}}{4}  \ell n\, \left(\frac{3+\sqrt{5}}{2} \right).
 \end{eqnarray}

\vskip.2cm
{\bf Concluding Remark :}
We conclude our present investigation by observing that several digamma functions for positive and negative fractional values have been deduced using our new identities \eqref{eq(dig-2.6)}, \eqref{eq(main-2.7)} and \eqref{eq(main-2.8)} in an analogous manner.


\begin{thebibliography}{99}

\bibitem{Al-saqabi}  Al-Saqabi, B.N. Kalla, S.L. and Srivastava, H.M.; A certain family of infinite series associated with Digamma functions, \textit{J. Math. Anal. Appl.}, {\bf 159} (1991), 361-372.

\bibitem{Alzer}  Alzer, H.; Sharp inequalities for digamma and polygamma functions, \textit{Forum Math.}, {\bf16} (2004), 181-221.

\bibitem {Batir}  Batir, N.; Some new inequalities for gamma and polygamma function, \textit{JIPAM. J. Inequal. Pure Appl. Math.}, {\bf 6(4)} (2005), Article 103, 9 p.

\bibitem{Boh}  B\"{o}hmer, E.; \textit{Differenzengleichungen und bestimmte integrale}, Leipzig, (1939).

\bibitem{Bor} Borwein, D. and Borwein, J.M.; On an intriguing integral and some series related to $\zeta(4)$,\textit{Proc. Amer. Math. Soc.}, {\bf 123} (1995), 1191-1198.

\bibitem {Clark} Clark, W.E. and Ismail, M.E.H.; Inequalities involving gamma and psi function, \textit{Anal. Appl.}, {\bf 1(129)} (2003), 129-140.

\bibitem{de}  De Doelder, P.J.; On some series containing $\psi(x)-\psi(y)$ and $\left(\psi(x)-\psi(y)\right)^{2}$ for certain values of $x$ and $y$, \textit{J. Comput. Appl. Math.}, {\bf 37} (1991), 125-141.

\bibitem{Erde}  Erd$\acute{e}$lyi, A. Magnus,  W. Oberhettinger, F. and Tricomi, F.G.; \textit{Higher Transcendental Functions, Vol.I} (Bateman Manuscript Project), McGraw-Hill, Book Co. Inc., New York, Toronto and London, (1953).

\bibitem{Gauss}  Gauss, C.F.; Disquisitiones generales circa seriem infinitam etc., \textit{Comm. Soc. reg. Sci. Gott. rec., Vol II}, (1813) pp. 1-46.; reprinted in Werke {\bf 3} (1866), 123-163.

\bibitem {Gosper}  Gosper, R.W.; $\int_{n/4}^{m/6}\log\Gamma(z)dz$, In special functions, q-series and related topics, \textit{Amer. Math. Soc.}, {\bf 14} (1997), 71-76.

\bibitem{Gradshteyn}  Gradshteyn, I.S. and  Ryzhik, I.M.; \textit{Table of integrals, series and products}, {\bf 8th ed.}, Academic Press Inc., San Diego, CA. 2014.

\bibitem{Grossman} Grossman, N.; Polygamma functions of arbitrary order, \textit{SIAM J. Math. Anal.}, {\bf 7} (1976), 366-372.

\bibitem{Jensen}  Jensen, J.L.W.V.; An elementary exposition of the theory of the Gamma function, \textit{Ann. Math.}, {\bf 17(3)} (1916), 124-166.

\bibitem{Lehmer} Lehmer, D.H.; Euler constants for arithmetical progressions, \textit{Acta Arith.}, {\bf 27} (1975), 125-142.

\bibitem{Lewin} Lewin, L.; \textit{Polygarithms and Associated Functions}, Elsevier, Amsterdam, 1981.

\bibitem{Mahler} Mahler, K.; Applications of a theorem of A. B. Shidlovski, \textit{Proc. Royal Soc. London Ser. A Math. Phys. Eng. Sci.}, {\bf 305} (1968), 149-173.

\bibitem{Murty} Murty, M.R. and Saradha, N.; Transcendental values of the digamma function, \textit{J. Number Theory}, {\bf 125} (2007), 298-318.

\bibitem{Nielsen} Nielsen, N.; \textit{Handbuch der theorie der gamma funktion}, Leipzig Druck und Verleg Von B.G. Teubner, 1906.

\bibitem{Sri1} Srivastava, H.M.; A simple algorithm for the evaluation of a class of generalized hypergeometric series, \textit{Stud. Appl. Math.}, {\bf 86} (1992), 79-86.

\bibitem{Sri2} Srivastava, H.M. and Choi, J.; \textit{Series Associated with the Zeta and Related Functions}, Kluwer, Dordrecht, 2001.

\bibitem{Qi1} Qi, F. and Chen, Ch.-P.; A complete monotonicity of the gamma function, \textit{RGMIA Res. Rep. Coll}, {\bf 7} (2007), Art. 1.

\bibitem{Qi2} Qi, F. and Chen, Ch.-P.; A complete monotonicity property of the gamma function, \textit{J. Math. Anal. Appl.}, {\bf 296} (2004), 603-607.

\bibitem{Qi3} Qi, F. and Guo, B.-P.; Complete monotonicities of functions involving the gamma and digamma functions, \textit{RGMIA Res. Rep. Coll}, {\bf 7} (2004), 63-72, Art. 8.

\bibitem{Qi4} Qi, F. Guo, B.-P. and  Chen, Ch.-P.; Some completely monotonic functions involving the gamma and polygamma functions, \textit{RGMIA Res. Rep. Coll}, {\bf 7} (2004), 31-36, Art. 5.

\bibitem{Qi}  Qiu, S.L. and Vuorinen, M.; Some properties of the gamma and Psi functions with applications, \textit{Math. Comp.}, {\bf 74} (2005), 723-742.

\bibitem{Qureshi} Qureshi, M.I. Jabee, Saima and Shadab, M.; Truncated Gauss hypergeometric series and its application in digamma function, \textit{(Communicated)}.

\bibitem{Wu}  Wu, T.-C. Leu, S.-H. Tu, S.-T. and Srivastava, H.M.; A certain class of infinite sums associated with Diagamma functions, \textit{Appl. Math. Comput.}, {\bf 105} (1999), 1-9.

\end{thebibliography}
\end{document}